\documentclass[12pt]{article}

\usepackage{amssymb}
\usepackage{luis}

\def\stackunder#1#2{\mathrel{\mathop{#2}\limits_{#1}}}
\def\QTR#1#2{{\csname#1\endcsname #2}}

\def\func#1{{\rm #1}}
\def\dsum{\mathop{\displaystyle \sum }}%
\def\limfunc#1{\mathop{\rm #1}}%

\newtheorem{theorem}{Theorem}[section]

\newtheorem{proposition}[theorem]{Proposition}
\newtheorem{definition}[theorem]{Definition}
\newtheorem{example}[theorem]{Example}

\newtheorem{remark}[theorem]{Remark}

\begin{document}

\title{On a relation between intrinsic and extrinsic Dirichlet forms for 
interacting particle systems}

\author{\textbf{Jos\'{e} L.~da Silva} \\
BiBoS, Bielefeld Univ., D-33615 Bielefeld, Germany \\            
CCM, Madeira Univ, P-9000 Funchal, Madeira Portugal\\
\and 
\textbf{Yuri G.~Kondratiev} \\
Inst.~Ang.~Math., Bonn Univ., D-53115 Bonn, Germany\\
BiBoS, Bielefeld Univ., D-33615 Bielefeld, Germany \\            
Inst.~Math., NASU, 252601 Kiev, Ukraine\\ 
\and 
\textbf{Michael R\"ockner} \\
Fakult\"at f\"ur Mathematik, Universit\"at Bielefeld,\\
D-33615 Bielefeld, Germany            
}
\maketitle

\begin{abstract}
In this paper we extend the result obtained in \cite{AKR97} (see also \cite
{AKR96}) on the representation of the intrinsic pre-Dirichlet form $\mathcal{%
E}_{\pi _{\sigma }}^{\Gamma }$ of the Poisson measure $\pi _{\sigma }$ in
terms of the extrinsic one $\mathcal{E}_{\pi _{\sigma },H_{\sigma }^{X}}^{P}$%
. More precisely, replacing $\pi _{\sigma }$ by a Gibbs measure $\mu $ on
the configuration space $\Gamma _{X}$ we derive a relation between the
intrinsic pre-Dirichlet form $\mathcal{E}_{\mu }^{\Gamma }$ of the measure $%
\mu $ and the extrinsic one $\mathcal{E}_{\mu ,H_{\sigma }^{X}}^{P}$. As a
consequence we prove the closability of $\mathcal{E}_{\mu }^{\Gamma }$ on $%
L^{2}(\Gamma _{X},\mu )$ under very general assumptions on the interaction
potential of the Gibbs measures $\mu $.
\end{abstract}

\section{Introduction}

In the recent papers \cite{AKR96}, \cite{AKM96}, \cite{AKR97}, and \cite
{AKR97a} analysis and geometry on configuration spaces $\Gamma _{X}$ over a
Riemannian manifold $X$, i.e, 
\[
\Gamma _{X}:=\{\gamma \subset X\,|\,|\gamma \cap K|<\infty \;\mathrm{%
for\;any\;compact\;}K\subset X\}, 
\]
was developed. One of the consequences of the discussed approach was a
description of the well-known equilibrium process on configuration spaces as
the Brownian motion associated with a Dirichlet form of the Poisson measure $%
\pi _{\sigma }$ with intensity measure $\sigma $ on $\mathcal{B}(X)$. This
form is canonically associated with the introduced geometry on configuration
spaces and is called \textit{intrinsic Dirichlet form} of the measure $\pi
_{\sigma }$.

On the other hand there is a well-known realization (canonical isomorphism)
of the Hilbert space $L^{2}(\Gamma _{X},\pi _{\sigma })$ and the
corresponding Fock space 
\[
\func{Exp}L^{2}(X,\sigma ):=\bigoplus_{n=0}^{\infty }\func{Exp}%
_{n}L^{2}(X,\sigma ),
\]
where $\func{Exp}_{n}L^{2}(X,\sigma )$ denotes the $n$-fold symmetric tensor
product of $L^{2}(X,\sigma )$ and $\func{Exp}_{0}L^{2}(X,\sigma ):=%
\QTR{mathbb}{C}$. This isomorphism produces natural operations in $%
L^{2}(\Gamma _{X},\pi _{\sigma })$ as images of the standard Fock space
operators, see e.g., \cite{KSSU97} and references therein. In particular, we
can consider the image of the annihilation operator from the Fock space as a
natural version of a ``gradient'' operator in $L^{2}(\Gamma _{X},\pi
_{\sigma })$. The differentiable structure in $L^{2}(\Gamma _{X},\pi
_{\sigma })$ which appears in this way we consider as \textit{external}
because it is produced via transportation from the Fock space.

As was shown in \cite[Section 5]{AKR97} the intrinsic Dirichlet form of the
measure $\pi _{\sigma }$ can be represented also in terms of the external
Dirichlet form $\mathcal{E}_{\pi _{\sigma },H_{\sigma }^{X}}^{P}$ with
coefficient $H_{\sigma }^{X}$ (the Dirichlet operator associated with $%
\sigma $ on $X$) which uses this external differentiable structure, i.e., 
\[
\int_{\Gamma }\langle \nabla ^{\Gamma }F(\gamma ),\nabla ^{\Gamma }G(\gamma
)\rangle _{T_{\gamma }\Gamma }d\pi _{\sigma }(\gamma )=\int_{\Gamma }(\nabla
^{P}F(\gamma ),H_{\sigma }^{X}\nabla ^{P}G(\gamma ))_{L^{2}(X,\sigma )}d\pi
_{\sigma }(\gamma ). 
\]
As a result we have a full spectral description of the corresponding
Dirichlet operator $H_{\pi _{\sigma }}^{\Gamma }$ which is the generator of
the equilibrium process on $\Gamma _{X}$.

If we change the Poisson measure $\pi _{\sigma }$ to a Gibbs measure $\mu $
on the configuration space $\Gamma _{X}$ which describes the equilibrium of
interacting particle systems, the corresponding intrinsic Dirichlet form can
still be used for constructing the corresponding stochastic dynamics (cf.~%
\cite[Section 5]{AKR97a}) or for constructing a quantum infinite particle
Hamiltonian in models of quantum fields theory, see \cite{AKR97b}.

The aim of this paper is to show that even for the interacting case there is
a transparent relation between the intrinsic Dirichlet form and the
extrinsic one, see Theorem \ref{5eq2}. The proof is based on the
Nguyen-Zessin characterization of Gibbs measure (cf.~\cite[Theorem 2]{NZ79}
or Proposition \ref{5eq48} below) which on a heuristic level can be
considered as a consequence of the Mecke identity (cf.~\cite[Satz 3.1]{Me67}%
), see Remark \ref{5eq49} below for more details.

As a consequence of the mentioned relation we prove the closability of the
pre-Dirichlet form $(\mathcal{E}_{\mu }^{\Gamma },\mathcal{F}C_{b}^{\infty }(%
\mathcal{D},\Gamma ))$ on $L^{2}(\Gamma _{X},\mu )$, where $\mu $ is a
tempered grand canonical Gibbs measure, see Section \ref{5eq6} for this
notion. It turns out that this result is obtained as a ``lifting'' of the
closable Dirichlet forms on $X$. We would like to emphasize that we achieve
this result under a general condition (see (\ref{5eq51}) below) on the
potential $\Phi $ which is not covered by condition (A.6) in \cite{O96}. The
closability is crucial (for physical reasons, see \cite{AKR97b}, and) for
applying the general theory of Dirichlet forms including the construction of
a corresponding diffusion process (cf.~\cite{MR92}) which models an infinite
particle system with (possibly) very singular interactions (cf.~\cite{AKR97a}%
).

Another motivation for deriving Theorem \ref{5eq2} is to use this result for
studying spectral properties of Hamiltonians of intrinsic Dirichlet forms
associated with Gibbs measures. This will be implemented in a forthcoming
paper.

\section{Preliminaries and Framework\label{5eq6}}

In this section we describe some facts about probability measures on
configuration spaces which are necessary later on.

Let $X$ be a connected, oriented $C^{\infty }$ (non-compact) Riemannian
manifold. For each point $x\in X$, the tangent space to $X$ at $x$ will be
denoted by $T_{x}X$; and the tangent bundle will be denoted by $TX=\cup
_{x\in X}T_{x}X$. The Riemannian metric on $X$ associates to each point $%
x\in X$ an inner product on $T_{x}X$ which we denote by $\left\langle \cdot
,\cdot \right\rangle _{T_{x}X}$ and the associated norm will be denoted by $%
|\cdot |_{T_{x}X}$. Let $m$ denote the volume element.

$\mathcal{O}(X)$ is defined as the family of all open subsets of $X$ and $%
\mathcal{B}(X)$ denotes the corresponding Borel $\sigma $-algebra. $\mathcal{%
O}_{c}(X)$ and $\mathcal{B}_{c}(X)$ denote the systems of all elements in $%
\mathcal{O}(X)$, $\mathcal{B}(X)$ respectively, which have compact closures.

Let $\Gamma :=\Gamma _{X}$ be the set of all locally finite subsets in $X$: 
\[
\Gamma _{X}:=\{\gamma \subset X\,|\,|\gamma \cap K|<\infty \;\mathrm{%
for\;any\;compact\;}K\subset X\}. 
\]

We will identify $\gamma $ with the positive integer-valued measure $%
\sum_{x\in \gamma }\varepsilon _{x}$. Then for any $\varphi \in C_{0}(X)$ we
have a functional $\Gamma \ni \gamma \mapsto \langle \varphi ,\gamma \rangle
=\sum_{x\in \gamma }\varphi (x)\in \Bbb{R}$. Here $C_{0}(X)$ is the set of
all real-valued continuous functions on $X$ with compact support. The space $%
\Gamma $ is endowed with the vague topology. Let $\mathcal{B}(\Gamma )$
denote the corresponding Borel $\sigma $-algebra. For $\Lambda \subset X$ we
sometimes use the shorthand $\gamma _{\Lambda }$ for $\gamma \cap \Lambda $.

For any $B\in \mathcal{B}(X)$ we define, as usual, $\Gamma \ni \gamma
\mapsto N_{B}(\gamma ):=\gamma (B)\in \Bbb{Z}_{+}\cup \{+\infty \}$. Then $%
\mathcal{B}(\Gamma )=\sigma (\{N_{\Lambda }|\Lambda \in \mathcal{O}%
_{c}(X)\}).$ For any $A\in \mathcal{B}(X)$ we also define $\mathcal{B}%
_{A}(\Gamma ):=\sigma (\{N_{B}|B\in \mathcal{B}_{c}(X),B\subset A\}).$

Let $d\sigma (x)=\rho (x)dm(x)$, where $\rho >0$ $m$-a.e.~such that $\rho ^{%
\frac{1}{2}}\in H_{loc}^{1,2}(X)$ (the Sobolev space of order 1 in $%
L^{2}(X,m)$) and $\rho \notin L^{1}(X,m)$. We recall that the Poisson
measure $\pi _{\sigma }$ (with intensity measure $\sigma $) on $(\Gamma ,%
\mathcal{B}(\Gamma ))$ is defined via its Laplace transform by 
\begin{equation}
\int_{\Gamma }e^{\langle \gamma ,\varphi \rangle }d\pi _{\sigma }(\gamma
)=\exp \left( \int_{X}(e^{\varphi (x)}-1)d\sigma (x)\right) ,\;\varphi \in
C_{0}(X),  \label{5eq18}
\end{equation}
see e.g.~\cite{AKR97}, \cite{GV68}, and \cite{S94}. Let us mention that if $%
\rho \in L^{1}(X,m)$, then we have a finite intensity measure $\sigma $ on $X
$, and in this case the corresponding measure $\pi _{\sigma }$ will be
concentrated on finite configurations. The latter can be considered as a
degenerated case which can be reduced to finite dimensional analysis on
every subset of $n$-particle configurations.

Let us briefly recall the definition of grand canonical Gibbs measures on $%
(\Gamma ,\mathcal{B}(\Gamma ))$. We adopt the notation in \cite{AKR97a}, and
refer the interested reader to the beautiful work by C.~Preston, \cite{P79},
but also \cite{P76}, and \cite{G79}.

A function $\Phi :\Gamma \rightarrow \Bbb{R}\cup \{+\infty \}$ will be
called a \textit{potential} iff for all $\Lambda \in \mathcal{B}_{c}(X)$ we
have $\Phi (\emptyset )=0$, $\Phi =1\!\!1_{\{N_{X}<\infty \}}\Phi $, and $%
\gamma \mapsto \Phi (\gamma _{\Lambda })$ is $\mathcal{B}_{\Lambda }(\Gamma
) $-measurable.

For $\Lambda \in \mathcal{B}_{c}(X)$ the \textit{conditional energy} $%
E_{\Lambda }^{\Phi }:\Gamma \rightarrow \Bbb{R}\cup \{+\infty \}$ is defined
by 
\begin{equation}
E_{\Lambda }^{\Phi }(\gamma ):=\left\{ 
\begin{array}{cl}
\dsum\limits_{\gamma ^{\prime }\subset \gamma ,\gamma ^{\prime }(\Lambda
)>0}\Phi (\gamma ^{\prime }) & \mathrm{if\,}\dsum\limits_{\gamma ^{\prime
}\subset \gamma ,\gamma ^{\prime }(\Lambda )>0}|\Phi (\gamma ^{\prime
})|<\infty , \\ 
&  \\ 
+\infty & \mathrm{otherwise,}
\end{array}
\right.  \label{5eq19}
\end{equation}
where the sum of the empty set is defined to be zero.

Later on we will use conditional energies which satisfy an additional
assumption, namely, the \textit{stability condition}, i.e., there exists $%
B\geq 0$ such that for any $\Lambda \in \mathcal{B}_{c}(X)$ and for all $%
\gamma \in \Gamma _{\Lambda }$%
\[
E_{\Lambda }^{\Phi }(\gamma )\geq -B|\gamma |. 
\]

\begin{definition}
For any $\Lambda \in \mathcal{O}_{c}(X)$ define for $\gamma \in \Gamma $ the
measure $\Pi _{\Lambda }^{\Phi }(\gamma ,\cdot )$ by 
\begin{eqnarray}
\Pi _{\Lambda }^{\sigma ,\Phi }(\gamma ,\Delta ) &\mbox{$:=$}%
&1\!\!1_{\{Z_{\Lambda }^{\sigma ,\Phi }<\infty \}}(\gamma )[Z_{\Lambda
}^{\sigma ,\Phi }(\gamma )]^{-1}\int_{\Gamma }1\!\!1_{\Delta }(\gamma
_{X\backslash \Lambda }+\gamma _{\Lambda }^{\prime })  \label{5eq21} \\
&&\cdot \exp [-E_{\Lambda }^{\Phi }(\gamma _{X\backslash \Lambda }+\gamma
_{\Lambda }^{\prime })]d\pi _{\sigma }(\gamma ^{\prime }),\;\Delta \in 
\mathcal{B}(\Gamma ),  \nonumber
\end{eqnarray}
where 
\begin{equation}
Z_{\Lambda }^{\sigma ,\Phi }(\gamma ):=\int_{\Gamma }\exp [-E_{\Lambda
}^{\Phi }(\gamma _{X\backslash \Lambda }+\gamma _{\Lambda }^{\prime })]d\pi
_{\sigma }(\gamma ^{\prime }).  \label{5eq22}
\end{equation}

A probability measure $\mu $ on $(\Gamma ,\mathcal{B}(\Gamma ))$ is called
grand canonical Gibbs measure with interaction potential $\Phi $ if for all $%
\Lambda \in \mathcal{O}_{c}(X)$%
\begin{equation}
\mu \Pi _{\Lambda }^{\Phi }=\mu .  \label{5eq23}
\end{equation}
Let $\mathcal{G}_{gc}(\sigma ,\Phi )$ denote the set of all such probability
measures $\mu $.
\end{definition}

\begin{remark}
\begin{enumerate}
\item  It is well-known that $(\Pi _{\Lambda }^{\sigma ,\Phi })_{\Lambda \in 
\mathcal{O}_{c}(X)}$ is a $(\mathcal{B}_{X\backslash \Lambda }(\Gamma
))_{\Lambda \in \mathcal{O}_{c}(X)}$-specification in the sense of 
\cite[Section 6]{P76} or \cite{P79}.

\item  For any $\gamma \in \Gamma $ the measure $\mu \Pi _{\Lambda }^{\Phi }$
in (\ref{5eq23}) is defined by 
\begin{equation}
(\mu \Pi _{\Lambda }^{\sigma ,\Phi })(\Delta ):=\int_{\Gamma }d\mu (\gamma
)\Pi _{\Lambda }^{\sigma ,\Phi }(\gamma ,\Delta ),\;\Delta \in \mathcal{B}%
(\Gamma )  \label{5eq24}
\end{equation}
and (\ref{5eq23}) are called Dobrushin-Landford-Ruelle (DLR) equations.
\end{enumerate}
\end{remark}

\section{Intrinsic geometry on Poisson space\label{5eq10}}

We recall some results to be used below from \cite{AKR97}, \cite{AKM96} to
which we refer for the corresponding proofs and more details.

A homeomorphism $\psi :X\rightarrow X$ defines a transformation of $\Gamma $
by 
\[
\Gamma \ni \gamma \mapsto \psi (\gamma )=\{\psi (x)|x\in \gamma
\}=\sum_{x\in \gamma }\varepsilon _{\psi (x)}. 
\]
Any vector field $v\in V_{0}(X)$ (i.e., the set of all smooth vector fields
on $X$ with compact support) defines (via the exponential mapping) a
one-parameter group $\psi _{t}^{v}$, $t\in \Bbb{R}$, of diffeomorphisms of $%
X $.

\begin{definition}
For $F:\Gamma \rightarrow \Bbb{R}$ we define the directional derivative
along the vector field $v$ as (provided the right hand side exists) 
\[
(\nabla _{v}^{\Gamma }F)(\gamma ):=\frac{d}{dt}F(\psi _{t}^{v}(\gamma
))|_{t=0}.
\]
\end{definition}

This definition applies to $F$ in the following class $\mathcal{F}%
C_{b}^{\infty }(\mathcal{D},\Gamma )$ of so-called smooth cylinder
functions. Let $\mathcal{D}:=C_{0}^{\infty }(X)$ (the set of all smooth
functions on $X$ with compact support). We define $\mathcal{F}C_{b}^{\infty
}(\mathcal{D},\Gamma )$ as the set of all functions on $\Gamma $ of the form 
\begin{equation}
F(\gamma )=g_{F}(\langle \gamma ,\varphi _{1}\rangle ,\ldots ,\langle \gamma
,\varphi _{N}\rangle ),\;\gamma \in \Gamma ,  \label{5eq20}
\end{equation}
where $\varphi _{1},\ldots ,\varphi _{N}\in \mathcal{D}$ and $g_{F}$ is from 
$C_{b}^{\infty }(\QTR{mathbb}{R}^{N})$. Clearly, $\mathcal{F}C_{b}^{\infty }(%
\mathcal{D},\Gamma )$ is dense in $L^{2}(\pi _{\sigma }):=L^{2}(\Gamma ,\pi
_{\sigma })$. For any $F\in \mathcal{F}C_{b}^{\infty }(\mathcal{D},\Gamma )$
we have 
\begin{equation}
(\nabla _{v}^{\Gamma }F)(\gamma )=\sum_{i=1}^{N}\frac{\partial g_{F}}{%
\partial s_{i}}(\langle \gamma ,\varphi _{1}\rangle ,\ldots ,\langle \gamma
,\varphi _{N}\rangle )\langle \gamma ,\nabla _{v}^{X}\varphi _{i}\rangle ,
\label{5eq16}
\end{equation}
where $x\mapsto (\nabla _{v}^{X}\varphi )(x)=\langle \nabla ^{X}\varphi
(x),v(x)\rangle _{TX}$ is the usual directional derivative on $X$ along the
vector field $v$ and $\nabla ^{X}$ denotes the gradient on $X$.

The logarithmic derivative of the measure $\sigma $ is given by the vector
field $\beta ^{\sigma }:=\nabla ^{X}\log \rho =\nabla ^{X}\rho /\rho $
(where $\beta ^{\sigma }=0$ on $\{\rho =0\}$). Then the logarithmic
derivative of $\sigma $ along $v$ is the function $x\mapsto \beta
_{v}^{\sigma }(x)=\langle \beta ^{\sigma }(x),v(x)\rangle _{T_{x}X}+\func{div%
}^{X}v(x)$, where $\limfunc{div}^{X}$ denotes the divergence on $X$
w.r.t.~the volume element $m$. Analogously, we define $\limfunc{div}_{\sigma
}^{X}$ as the divergence on $X$ w.r.t.~$\sigma $, i.e., $\limfunc{div}%
_{\sigma }^{X}$ is the dual operator on $L^{2}(\sigma ):=L^{2}(X,\sigma )$
of $\nabla ^{X}$.

\begin{definition}
For any $v\in V_{0}(X)$ we define the logarithmic derivative of $\pi
_{\sigma }$ along $v$ as the following function on $\Gamma :$%
\begin{equation}
\Gamma \ni \gamma \mapsto B_{v}^{\pi _{\sigma }}(\gamma ):=\langle \beta
_{v}^{\sigma },\gamma \rangle =\int_{X}[\langle \beta ^{\sigma
}(x),v(x)\rangle _{T_{x}X}+\mathrm{div}^{X}v(x)]d\gamma (x).  \label{5eq25}
\end{equation}
\end{definition}

\begin{theorem}
For all $F,G\in $ $\mathcal{F}C_{b}^{\infty }(\mathcal{D},\Gamma )$ and any $%
v\in V_{0}(X)$ the following integration by parts formula for $\pi _{\sigma }
$ holds: 
\begin{equation}
\int_{\Gamma }\nabla _{v}^{\Gamma }FGd\pi _{\sigma }=-\int_{\Gamma }F\nabla
_{v}^{\Gamma }Gd\pi _{\sigma }-\int_{\Gamma }FGB_{v}^{\pi _{\sigma }}d\pi
_{\sigma },  \label{5eq26}
\end{equation}
or $(\nabla _{v}^{\Gamma })^{*}=-\nabla _{v}^{\Gamma }-B_{v}^{\pi _{\sigma }}
$, as an operator equality on the domain $\mathcal{F}C_{b}^{\infty }(%
\mathcal{D},\Gamma )$ in $L^{2}(\pi _{\sigma })$.
\end{theorem}

\begin{definition}
We introduce the tangent space $T_{\gamma }\Gamma $ to the configuration
space $\Gamma $ at the point $\gamma \in \Gamma $ as the Hilbert space of $%
\gamma $-square-integrable sections (measurable vector fields) $%
V:X\rightarrow TX$ with scalar product $\langle V^{1},V^{2}\rangle
_{T_{\gamma }\Gamma }=\int_{X}\langle V^{1}(x),V^{2}(x)\rangle
_{T_{x}X}d\gamma (x)$, $V^{1},V^{2}\in T_{\gamma }\Gamma =L^{2}(X\rightarrow
TX;\gamma )$. The corresponding tangent bundle is denoted by $T\Gamma $.
\end{definition}

The intrinsic gradient of a function $F\in \mathcal{F}C_{b}^{\infty }(%
\mathcal{D},\Gamma )$ is a mapping $\Gamma \ni \gamma \mapsto (\nabla
^{\Gamma }F)(\gamma )\in T_{\gamma }\Gamma $ such that $(\nabla _{v}^{\Gamma
}F)(\gamma )=\langle \nabla ^{\Gamma }F(\gamma ),v\rangle _{T_{\gamma
}\Gamma }$ for any $v\in V_{0}(X)$. Furthermore, by (\ref{5eq16}), if $F$ is
given by (\ref{5eq20}), we have for $\gamma \in \Gamma $, $x\in X$%
\begin{equation}
(\nabla ^{\Gamma }F)(\gamma ;x)=\sum_{i=1}^{N}\frac{\partial g_{F}}{\partial
s_{i}}(\langle \varphi _{1},\gamma \rangle ,\ldots ,\langle \varphi
_{N},\gamma \rangle )\nabla ^{X}\varphi _{i}(x).  \label{5eq44}
\end{equation}

\begin{definition}
For a measurable vector field $V:\Gamma \rightarrow T\Gamma $ the divergence 
$\limfunc{div}_{\pi _{\sigma }}^{\Gamma }V$ is defined via the duality
relation for all $F\in \mathcal{F}C_{b}^{\infty }(\mathcal{D},\Gamma )$ by 
\begin{equation}
\int_{\Gamma }\langle V_{\gamma },\nabla ^{\Gamma }F(\gamma )\rangle
_{T_{\gamma }\Gamma }d\pi _{\sigma }(\gamma )=-\int_{\Gamma }F(\gamma )(%
\mathrm{div}_{\pi _{\sigma }}^{\Gamma }V)(\gamma )d\pi _{\sigma }(\gamma ),
\label{5eq27}
\end{equation}
provided it exists (i.e., provided 
\[
F\mapsto \int_{\Gamma }\langle V_{\gamma },\nabla ^{\Gamma }F(\gamma
)\rangle _{T_{\gamma }\Gamma }d\pi _{\sigma }(\gamma ) 
\]
is continuous on $L^{2}(\pi _{\sigma })$).
\end{definition}

\begin{proposition}
\label{5eq33}For any vector field $V=G\,v$, where $G\in \mathcal{F}%
C_{b}^{\infty }(\mathcal{D},\Gamma )$, $v\in V_{0}(X)$ we have 
\begin{equation}
(\mathrm{div}_{\pi _{\sigma }}^{\Gamma }V)(\gamma )=\langle (\nabla ^{\Gamma
}G)(\gamma ),v\rangle _{T_{\gamma }\Gamma }+G(\gamma )B_{v}^{\pi _{\sigma
}}(\gamma ).  \label{5eq28}
\end{equation}
\end{proposition}

For any $F,G\in \mathcal{F}C_{b}^{\infty }(\mathcal{D},\Gamma )$ we
introduce the \textit{pre-Dirichlet form} which is generated by the
intrinsic gradient $\nabla ^{\Gamma }$ as 
\begin{equation}
\mathcal{E}_{\pi _{\sigma }}^{\Gamma }(F,G)=\int_{\Gamma }\langle (\nabla
^{\Gamma }F)(\gamma ),(\nabla ^{\Gamma }G)(\gamma )\rangle _{T_{\gamma
}\Gamma }d\pi _{\sigma }(\gamma ).  \label{5eq29}
\end{equation}

We will also need the \textit{classical pre-Dirichlet} form for the
intensity measure $\sigma $ which is given as $\mathcal{E}_{\sigma
}^{X}(\varphi ,\psi )=\int_{X}\langle \nabla ^{X}\varphi ,\nabla ^{X}\psi
\rangle _{TX}d\sigma $ for any $\varphi ,\psi \in \mathcal{D}$. This form is
associated with the \textit{Dirichlet operator} $H_{\sigma }^{X}$ which is
given on $\mathcal{D}$ by $H_{\sigma }^{X}\varphi (x):=-\triangle
^{X}\varphi (x)-\langle \beta ^{\sigma }(x),\nabla ^{X}\varphi (x)\rangle
_{T_{x}X}$ and which satisfies $\mathcal{E}_{\sigma }^{X}(\varphi ,\psi
)=(H_{\sigma }^{X}\varphi ,\psi )_{L^{2}(\sigma )}$, $\varphi ,\psi \in 
\mathcal{D}$, see e.g.~\cite{BK88} and \cite{MR92}.

For any $F\in \mathcal{F}C_{b}^{\infty }(\mathcal{D},\Gamma )$, $(\nabla
^{\Gamma }\nabla ^{\Gamma }F)(\gamma ,x,y)\in T_{\gamma }\Gamma \otimes
T_{\gamma }\Gamma $ and we can define the $\Gamma $-\textit{Laplacian} $%
(\triangle ^{\Gamma }F)(\gamma ):=\func{Tr}(\nabla ^{\Gamma }\nabla ^{\Gamma
}F)(\gamma )\in \mathcal{F}C_{b}^{\infty }(\mathcal{D},\Gamma )$. We
introduce a differential operator in $L^{2}(\pi _{\sigma })$ on the domain $%
\mathcal{F}C_{b}^{\infty }(\mathcal{D},\Gamma )$ by the formula 
\begin{equation}
(H_{\pi _{\sigma }}^{\Gamma }F)(\gamma )=-\triangle ^{\Gamma }F(\gamma
)-\langle \mathrm{div}_{\sigma }^{X}(\nabla ^{\Gamma }F)(\gamma ,\cdot
),\gamma \rangle .  \label{5eq30}
\end{equation}

\begin{theorem}
\label{5eq35}The operator $H_{\pi _{\sigma }}^{\Gamma }$ is associated with
the intrinsic Dirichlet form $\mathcal{E}_{\pi _{\sigma }}^{\Gamma }$, i.e., 
\begin{equation}
\mathcal{E}_{\pi _{\sigma }}^{\Gamma }(F,G)=(H_{\pi _{\sigma }}^{\Gamma
}F,G)_{L^{2}(\pi _{\sigma })},  \label{5eq31}
\end{equation}
or $H_{\pi _{\sigma }}^{\Gamma }=-\limfunc{div}_{\pi _{\sigma }}^{\Gamma
}\nabla ^{\Gamma }$ on $\mathcal{F}C_{b}^{\infty }(\mathcal{D},\Gamma )$. We
call $H_{\pi _{\sigma }}^{\Gamma }$ the intrinsic Dirichlet operator of the
measure $\pi _{\sigma }$.
\end{theorem}

\section{Extrinsic geometry on Poisson space\label{5eq13}}

We recall the extrinsic geometry on $L^{2}(\pi _{\sigma })$ based on the
isomorphism with the Fock space. Our approach is based on \cite{KSS96} but
we should also mention \cite{BLL95}, \cite{IK88}, \cite{KSSU97}, \cite{NV95}%
, \cite{P95} for related considerations and references therein. For proofs
of the results stated below in this section, we refer to \cite[Sect.~5]
{AKR97}.

Let us define another ``gradient'' on functions $F:\Gamma \rightarrow 
\QTR{mathbb}{R}$ This gradient $\nabla ^{P}$ has specific useful properties
on Poissonian spaces. We will call $\nabla ^{P}$ the \textit{Poissonian
gradient}. To this end we consider as the tangent space to $\Gamma $ at any
point $\gamma \in \Gamma $ the same space $L^{2}(\sigma )$ and define a
mapping $\mathcal{F}C_{b}^{\infty }(\mathcal{D},\Gamma )\ni F\mapsto \nabla
^{P}F\in L^{2}(\pi _{\sigma })\otimes L^{2}(\sigma )$ by 
\begin{equation}
(\nabla ^{P}F)(\gamma ,x):=F(\gamma +\varepsilon _{x})-F(\gamma ),\;\gamma
\in \Gamma ,x\in X.  \label{5eq17}
\end{equation}
We stress that the transformation $\Gamma \ni \gamma \mapsto \gamma
+\varepsilon _{x}\in \Gamma $ is $\pi _{\sigma }$-a.e.~well-defined because $%
\pi _{\sigma }(\{\gamma \in \Gamma |x\in \gamma \})=0$ for any $x\in X$. The
directional derivative is then defined as 
\begin{equation}
(\nabla _{\varphi }^{P}F)(\gamma )=(\nabla ^{P}F(\gamma ,\cdot ),\varphi
)_{L^{2}(\sigma )}=\int_{X}[F(\gamma +\varepsilon _{x})-F(\gamma )]\varphi
(x)d\sigma (x),\;\varphi \in \mathcal{D}.  \label{5eq36}
\end{equation}
The Poissonian gradient $\nabla ^{P}$ yields (via a corresponding
``integration by parts'' formula) an orthogonal system of Charlier
polynomials on the Poisson space $(\Gamma ,\mathcal{B}(\Gamma ),\pi _{\sigma
})$.

For any $n\in \QTR{mathbb}{N}$ and all $\varphi \in \mathcal{D}$ we
introduce a function in $L^{2}(\pi _{\sigma })$ by 
\begin{equation}
Q_{n}^{\pi _{\sigma }}(\gamma ;\varphi ^{\otimes n}):=((\nabla _{\varphi
}^{P})^{*n}1)(\gamma ),  \label{5eq37}
\end{equation}
and define $Q_{0}^{\pi _{\sigma }}:=1$. Due to the kernel theorem 
\cite[Chap.~1]{BK88} these functions have the representation $Q_{n}^{\pi
_{\sigma }}(\gamma ;\varphi ^{\otimes n})=\langle Q_{n}^{\pi _{\sigma
}}(\gamma ),\varphi ^{\otimes n}\rangle $, with gene\-ralized symmetric
kernels $\Gamma \ni \gamma \mapsto Q_{n}^{\pi _{\sigma }}(\gamma )\in \func{%
Exp}_{n}\mathcal{D}^{\prime }$, $n\in \Bbb{N}$. Here and below by $\func{Exp}%
_{n}E$ we denote the $n$-th symmetric tensor power of a linear space $E$.
Then for any smooth kernel $\varphi ^{(n)}\in \func{Exp}_{n}\mathcal{D}%
^{\otimes n}$ we introduce the function $Q_{n}^{\pi _{\sigma }}(\gamma
;\varphi ^{(n)}):=\langle Q_{n}^{\pi _{\sigma }}(\gamma ),\varphi
^{(n)}\rangle $ such that for all $\varphi ^{(n)}\in \func{Exp}_{n}\mathcal{D%
}^{\otimes n}$, $\psi ^{(m)}\in \func{Exp}_{m}\mathcal{D}^{\otimes m}$%
\begin{equation}
\int_{\Gamma }Q_{n}^{\pi _{\sigma }}(\gamma ;\varphi ^{(n)})Q_{m}^{\pi
_{\sigma }}(\gamma ;\psi ^{(m)})d\pi _{\sigma }(\gamma )=\delta
_{nm}n!(\varphi ^{(n)},\psi ^{(m)})_{L^{2}(\sigma ^{\otimes n})}.
\label{5eq38}
\end{equation}
Hence (\ref{5eq37}) extends to the case of kernels from the so-called $n$%
-particle Fock space $\func{Exp}_{n}L^{2}(\sigma )$, $n\in \Bbb{N}$, and we
set $\func{Exp}_{0}L^{2}(\sigma ):=\Bbb{R}$.

As usual the symmetric Fock space over the Hilbert space $L^{2}(\sigma )$ is
defined as $\func{Exp}L^{2}(\sigma ):=\oplus _{n=0}^{\infty }\func{Exp}%
_{n}L^{2}(\sigma )$, see e.g.~\cite{BK88} and \cite{HKPS93}. The square of
the norm of a vector $(f^{(n)})_{n=0}^{\infty }\in \func{Exp}L^{2}(\sigma )$
is given by $\sum_{n=0}^{\infty }n!\left\| f^{(n)}\right\| _{L^{2}(\sigma
^{\otimes n})}^{2}$. For any $F\in L^{2}(\pi _{\sigma })$ there exists such
a Fock vector, so that we have the following chaos decomposition 
\begin{equation}
F(\gamma )=\sum_{n=0}^{\infty }Q_{n}^{\pi _{\sigma }}(\gamma ;f^{(n)}),
\label{5eq39}
\end{equation}
and moreover $\left\| F\right\| _{L^{2}(\pi _{\sigma
})}^{2}=\sum_{n=0}^{\infty }n!\left\| f^{(n)}\right\| _{L^{2}(\sigma
^{\otimes n})}^{2.}$. And vice versa, by (\ref{5eq39}) any Fock vector
generates a function from $L^{2}(\pi _{\sigma })$. This produces an
isomorphism between $L^{2}(\pi _{\sigma })$ and $\func{Exp}L^{2}(\sigma )$.

\begin{remark}
In probability theory the functions $Q_{n}^{\pi _{\sigma }}(\gamma ;f^{(n)})$
are called the $n$-multiple stochastic integrals of $f^{(n)}$ with respect
to the compensated Poisson process generated by the Poisson measure $\pi
_{\sigma }$, see e.g.~\cite{NV95}.

There is an alternative approach to the chaos decomposition on the Poisson
space which uses the concept of generalized Appell systems, see e.g.~\cite
{KSS96}.
\end{remark}

The following proposition shows that the operators $\nabla _{\varphi }^{P}$
and $\nabla _{\varphi }^{P*}$ play the role of the annihilation
resp.~creation operators in the Fock space $\func{Exp}L^{2}(\sigma )$.

\begin{proposition}
For all $\varphi ,\psi \in \mathcal{D}$, $n\in \Bbb{N}$ the following
formulas hold 
\begin{equation}
\nabla _{\psi }^{P}Q_{n}^{\pi _{\sigma }}(\gamma ;\varphi ^{\otimes
n})=n(\varphi ,\psi )_{L^{2}(\sigma )}Q_{n-1}^{\pi _{\sigma }}(\gamma
;\varphi ^{\otimes (n-1)})  \label{5eq40}
\end{equation}
\begin{equation}
\nabla _{\psi }^{P*}Q_{n}^{\pi _{\sigma }}(\gamma ;\varphi ^{\otimes
n})=Q_{n+1}^{\pi _{\sigma }}(\gamma ;\varphi ^{\otimes n}\hat{\otimes}\psi
),\;\gamma \in \Gamma ,  \label{5eq15}
\end{equation}
where $\varphi ^{\otimes n}\hat{\otimes}\psi $ means the symmetric tensor
product of $\varphi ^{\otimes n}$ and $\psi $.
\end{proposition}

Next we give an explicit expression for the adjoint of the Poissonian
gradient $\nabla ^{P*}$.

\begin{proposition}
For any function $F\in L^{1}(\pi _{\sigma })\otimes L^{1}(\sigma )$ we have $%
F\in D(\nabla ^{P*})$ and the following equality holds 
\begin{equation}
(\nabla ^{P*}F)(\gamma )=\int_{X}F(\gamma -\varepsilon _{x},x)d\gamma
(x)-\int_{X}F(\gamma ,x)d\sigma (x),\;\gamma \in \Gamma ,  \label{5eq45}
\end{equation}
provided the right hand side of (\ref{5eq45}) is in $L^{2}(\pi _{\sigma })$.
\end{proposition}

\noindent \textbf{Proof.} For $X=\QTR{mathbb}{R}^{d}$ this proposition was
proved in \cite{KSSU97}. Let $G\in \mathcal{F}C_{b}^{\infty }(\mathcal{D}%
,\Gamma )$ be given. Then an application of (\ref{5eq17}) gives 
\begin{eqnarray}
(\nabla ^{P}G,F)_{L^{2}(\pi _{\sigma })\otimes L^{2}(\sigma )}
&=&\int_{X}\int_{\Gamma }G(\gamma +\varepsilon _{x})F(\gamma ,x)d\pi
_{\sigma }(\gamma )d\sigma (x)  \nonumber \\
&&-\int_{X}\int_{\Gamma }G(\gamma )F(\gamma ,x)d\pi _{\sigma }(\gamma
)d\sigma (x).  \label{5eq46}
\end{eqnarray}
Now we use the Mecke identity, see e.g., \cite[Satz 3.1]{Me67} 
\begin{equation}
\int_{\Gamma }\left( \int_{X}h(\gamma ,x)d\gamma (x)\right) d\pi _{\sigma
}(\gamma )=\int_{X}\int_{\Gamma }h(\gamma +\varepsilon _{x},x)d\pi _{\sigma
}(\gamma )d\sigma (x),  \label{5eq47}
\end{equation}
where $h$ is any non-negative, $\mathcal{B}(\Gamma )\times \mathcal{B}(X)$%
-measurable function. By (\ref{5eq47}) the right hand side of (\ref{5eq46})
transforms into 
\[
\int_{\Gamma }G(\gamma )\left[ \int_{X}F(\gamma -\varepsilon _{x},x)d\gamma
(x)-\int_{X}F(\gamma ,x)d\sigma (x)\right] d\pi _{\sigma }(\gamma ), 
\]
which proves the proposition.\hfill $\blacksquare \bigskip $

For any contraction $B$ in $L^{2}(\sigma )$ it is possible to define an
operator $\func{Exp}B$ as a contraction in $\func{Exp}L^{2}(\sigma )$ which
in any $n$-particle subspace $\func{Exp}_{n}L^{2}(\sigma )$ is given by $%
B\otimes \cdots \otimes B$ ($n$ times). For any positive self-adjoint
operator $A$ in $L^{2}(\sigma )$ (with $\mathcal{D}\subset D(A)$) we have a
contraction semigroup $e^{-tA}$, $t\geq 0$, hence it is possible to
introduce the second quantization operator $d\func{Exp}A$ as the generator
of the semigroup $\func{Exp}(e^{-tA})$, $t\geq 0$, i.e., $\func{Exp}%
(e^{-tA})=\exp (-td\func{Exp}A)$, see e.g., \cite{RS75}. We denote by $%
H_{A}^{P}$ the image of the operator $d\func{Exp}A$ in the Poisson space $%
L^{2}(\pi _{\sigma })$ under the described isomorphism.

\begin{proposition}
Let $\mathcal{D}\subset D(A)$. Then the symmetric bilinear form
corresponding to the operator $H_{A}^{P}$ has the following form 
\begin{equation}
(H_{A}^{P}F,G)_{L^{2}(\pi _{\sigma })}=\int_{\Gamma }(\nabla ^{P}F(\gamma
),A\nabla ^{P}G(\gamma ))_{L^{2}(\sigma )}d\pi _{\sigma }(\gamma ),\;F,G\in 
\mathcal{F}C_{b}^{\infty }(\mathcal{D},\Gamma ).  \label{5eq41}
\end{equation}
The right hand side of (\ref{5eq41}) is called the ``Poissonian
pre-Dirichlet form'' with coefficient operator $A$ and is denoted by $%
\mathcal{E}_{\pi _{\sigma },A}^{P}$.
\end{proposition}

Let us consider the special case of the second quantization operator $d\func{%
Exp}A$, where the one-particle operator $A$ coincides with the Dirichlet
operator $H_{\sigma }^{X}$ generated by the measure $\sigma $ on $X$. Then
we have the following theorem which relates the intrinsic Dirichlet operator 
$H_{\pi _{\sigma }}^{\Gamma }$ and the operator $H_{H_{\sigma }^{X}}^{P}$.

\begin{theorem}
\label{5eq42}$H_{\pi _{\sigma }}^{\Gamma }=H_{H_{\sigma }^{X}}^{P}$ on $%
\mathcal{F}C_{b}^{\infty }(\mathcal{D},\Gamma )$. In particular, for all $%
F,G\in \mathcal{F}C_{b}^{\infty }(\mathcal{D},\Gamma )$%
\begin{eqnarray}
&&\int_{\Gamma }\langle \nabla ^{\Gamma }F(\gamma ),\nabla ^{\Gamma
}G(\gamma )\rangle _{T_{\gamma }\Gamma }d\pi _{\sigma }(\gamma )  \nonumber
\\
&=&\int_{\Gamma }(\nabla ^{P}F(\gamma ),H_{\sigma }^{X}\nabla ^{P}G(\gamma
))_{L^{2}(\sigma )}d\pi _{\sigma }(\gamma )  \nonumber \\
&=&\int_{\Gamma }\int_{X}\langle \nabla ^{X}\nabla ^{P}F(\gamma ,x),\nabla
^{X}\nabla ^{P}G(\gamma ,x)\rangle _{T_{x}X}d\sigma (x)d\pi _{\sigma
}(\gamma ).  \label{5eq43}
\end{eqnarray}
\end{theorem}

\begin{remark}
Theorem \ref{5eq42} gives full information about the spectrum of the
intrinsic Dirichlet operator $H_{\pi _{\sigma }}^{\Gamma }$ on the Poisson
space in terms of the underlying Dirichlet operator $H_{\sigma }^{X}$ coming
from the intensity measure $\sigma $. In the following section we obtain an
analogue of Theorem \ref{5eq42} for the class of measures $\mathcal{G}%
_{gc}^{1}(\sigma ,\Phi )$ (see (\ref{5eq66}) below) which is the aim of this
paper.
\end{remark}

\section{Relation between intrinsic and extrinsic \protect\linebreak
Dirichlet forms\label{5eq1}}

Here we consider the class of measures $\mathcal{G}_{gc}^{1}(\sigma ,\Phi )$
consisting of all $\mu \in \mathcal{G}_{gc}(\sigma ,\Phi )$ such that 
\[
\int_{\Gamma }\gamma (K)d\mu (\gamma )<\infty \;\mathrm{for\;all\;compact\;}%
K\subset X. 
\]
We define for any $\mu \in \mathcal{G}_{gc}^{1}(\sigma ,\Phi )$ the
pre-Dirichlet form $\mathcal{E}_{\mu }^{\Gamma }$ by 
\begin{equation}
\mathcal{E}_{\mu }^{\Gamma }(F,G):=\int_{\Gamma }\langle \nabla ^{\Gamma
}F(\gamma ),\nabla ^{\Gamma }G(\gamma )\rangle _{T_{\gamma }\Gamma }d\mu
(\gamma ),\;F,G\in \mathcal{F}C_{b}^{\infty }(\mathcal{D},\Gamma ).
\label{5eq34}
\end{equation}

After all our preparations we are now going to prove an analogue of (\ref
{5eq43}) for $\mu \in \mathcal{G}_{gc}^{1}(\sigma ,\Phi )$. We would like to
emphasize that the corresponding formula (\ref{5eq3}) is not obtained from (%
\ref{5eq43}) by just replacing $\pi _{\sigma }$ by $\mu \in \mathcal{G}%
_{gc}^{1}(\sigma ,\Phi )$. The essential difference is, in addition, an
extra factor involving the conditional energy $E_{\Lambda }^{\Phi }$.

\begin{theorem}
\label{5eq2}For any $\mu \in \mathcal{G}_{gc}^{1}(\sigma ,\Phi )$, we have
for all $F,G\in \mathcal{F}C_{b}^{\infty }(\mathcal{D},\Gamma )$%
\begin{eqnarray}
\mathcal{E}_{\mu }^{\Gamma }(F,G) &=&\int_{\Gamma }\langle \nabla ^{\Gamma
}F(\gamma ),\nabla ^{\Gamma }G(\gamma )\rangle _{T_{\gamma }\Gamma }d\mu
(\gamma )  \label{5eq3} \\
&=&\int_{\Gamma }\int_{X}\langle \nabla ^{X}\nabla ^{P}F(\gamma ,x),\nabla
^{X}\nabla ^{P}G(\gamma ,x)\rangle _{T_{x}X}\,e^{-E_{\{x\}}^{\Phi }(\gamma
+\varepsilon _{x})}d\sigma (x)d\mu (\gamma )  \nonumber
\end{eqnarray}
\end{theorem}

\noindent \textbf{Proof.} Let us take any $F\in \mathcal{F}C_{b}^{\infty }(%
\mathcal{D},\Gamma )$ of the form (\ref{5eq20}). Then given $\gamma \in
\Gamma $ and $x\in X$ (\ref{5eq17}) implies that 
\begin{eqnarray*}
\nabla ^{X}\nabla ^{P}F(\gamma ,x) &=&\nabla ^{X}F(\gamma +\varepsilon _{x})
\\
&=&\sum_{i=1}^{N}\frac{\partial g_{F}}{\partial s_{i}}(\langle \varphi
_{1},\gamma \rangle +\varphi _{1}(x),\ldots ,\langle \varphi _{N},\gamma
\rangle +\varphi _{N}(x))\nabla ^{X}\varphi _{i}(x).
\end{eqnarray*}
Let us define $\hat{F}_{i}(\gamma ):=\frac{\partial g_{F}}{\partial s_{i}}%
(\langle \varphi _{1},\gamma \rangle ,\ldots ,\langle \varphi _{N},\gamma
\rangle )$, $i=1,\ldots ,N$. Obviously, it is enough to prove the equality (%
\ref{5eq3}) for $F=G$. Thus, inserting the result of $\nabla ^{X}\nabla
^{P}F(\gamma ,x)$ into the right hand side of (\ref{5eq3}) we obtain 
\begin{equation}
\int_{\Gamma }\int_{X}\sum_{i,j=1}^{N}\langle \nabla ^{X}\varphi
_{i}(x),\nabla ^{X}\varphi _{j}(x)\rangle _{T_{x}X}\hat{F}_{i}(\gamma
+\varepsilon _{x})\hat{F}_{j}(\gamma +\varepsilon _{x})e^{-E_{\{x\}}^{\Phi
}(\gamma +\varepsilon _{x})}d\sigma (x)d\mu (\gamma ).  \label{5eq5}
\end{equation}
Then we need the following useful proposition which generalizes the Mecke
identity to measures in $\mathcal{G}_{gc}(\sigma ,\Phi )$, see \cite{NZ79}, 
\cite{MMW79}.

\begin{proposition}
\label{5eq48}Let $h:\Gamma \times X\rightarrow \QTR{mathbb}{R}_{+}$ be $%
\mathcal{B}(\Gamma )\times \mathcal{B}(X)$-measurable, and let $\mu \in 
\mathcal{G}_{gc}(\sigma ,\Phi )$. Then we have 
\begin{equation}
\int_{\Gamma }\left( \int_{X}h(\gamma ,x)d\gamma (x)\right) d\mu (\gamma
)=\int_{X}\int_{\Gamma }h(\gamma +\varepsilon _{x},x)e^{-E_{\{x\}}^{\Phi
}(\gamma +\varepsilon _{x})}d\mu (\gamma )d\sigma (x).  \label{5eq4}
\end{equation}
\end{proposition}

Using this proposition we transform (\ref{5eq5}) into 
\[
\int_{\Gamma }\sum_{i,j=1}^{N}\hat{F}_{i}(\gamma )\hat{F}_{j}(\gamma
)\langle \langle \nabla ^{X}\varphi _{i}(\cdot ),\nabla ^{X}\varphi
_{j}(\cdot )\rangle _{TX},\gamma \rangle d\mu (\gamma ). 
\]
On the other hand using (\ref{5eq44}) we obtain 
\[
\langle \nabla ^{\Gamma }F(\gamma ),\nabla ^{\Gamma }G(\gamma )\rangle
_{T\Gamma }=\sum_{i,j=1}^{N}\hat{F}_{i}(\gamma )\hat{F}_{j}(\gamma )\langle
\langle \nabla ^{X}\varphi _{i}(\cdot ),\nabla ^{X}\varphi _{j}(\cdot
)\rangle _{TX},\gamma \rangle . 
\]
Therefore the equality on the dense $\mathcal{F}C_{b}^{\infty }(\mathcal{D}%
,\Gamma )$ is valid which proves the theorem.\hfill $\blacksquare $

\begin{remark}
\label{5eq49}For the reader's convenience let us give a heuristic proof of
the Nguyen-Zessin characterization of Gibbs measures in (\ref{5eq4}) which
really is a consequence of the Mecke identity (cf.~(\ref{5eq47})). Indeed,
let us write (heuristically) 
\[
d\mu (\gamma )=\frac{1}{Z^{\sigma ,\Phi }}e^{-E^{\Phi }(\gamma )}d\pi
_{\sigma }(\gamma ). 
\]
Then the function $E_{\{x\}}^{\Phi }(\gamma +\varepsilon _{x})=E^{\Phi
}(\gamma +\varepsilon _{x})-E^{\Phi }(\gamma )$ informally is the variation
of the potential energy $E^{\Phi }(\gamma )$ when we add to the
configuration $\gamma $ an additional point $x\in X$. Using this
representation we have 
\begin{eqnarray*}
&&\int_{X}\int_{\Gamma }h(\gamma +\varepsilon _{x},x)e^{-E_{\{x\}}^{\Phi
}(\gamma +\varepsilon _{x})}d\mu (\gamma )d\sigma (x) \\
&=&(Z^{\sigma ,\Phi })^{-1}\int_{X}\int_{\Gamma }h(\gamma +\varepsilon
_{x},x)e^{-E_{\{x\}}^{\Phi }(\gamma +\varepsilon _{x})}d\pi _{\sigma
}(\gamma )d\sigma (x).
\end{eqnarray*}
Then we use the Mecke identity to transform the right hand side of the above
equality into 
\[
(Z^{\sigma ,\Phi })^{-1}\int_{\Gamma }\left( \int_{X}h(\gamma ,x)e^{-E^{\Phi
}(\gamma )}d\gamma (x)\right) d\pi _{\sigma }(\gamma )=\int_{\Gamma }\left(
\int_{X}h(\gamma ,x)d\gamma (x)\right) d\mu (\gamma ). 
\]
The rigorous proof in \cite{NZ79} is obtained as a formalization of the
heuristic computations above.
\end{remark}

\section{Closability of intrinsic Dirichlet forms\label{5eq50}}

In this section we will prove the closability of the intrinsic Dirichlet
form $(\mathcal{E}_{\mu }^{\Gamma },\mathcal{F}C_{b}^{\infty }(\mathcal{D}%
,\Gamma ))$ on $L^{2}(\mu ):=L^{2}(\Gamma ,\mu )$ for all $\mu \in \mathcal{G%
}_{gc}^{1}(\sigma ,\Phi )$, using the integral representation (\ref{5eq3})
in Theorem \ref{5eq2}. The closability of $(\mathcal{E}_{\mu }^{\Gamma },%
\mathcal{F}C_{b}^{\infty }(\mathcal{D},\Gamma ))$ over $\Gamma $ is implied
by the closability of an appropriate family of pre-Dirichlet forms over $X$.
Let us describe this more precisely.

We define new intensity measures on $X$ by $d\sigma _{\gamma }(x):=\rho
_{\gamma }(x)dm(x)$, where 
\begin{equation}
\rho _{\gamma }(x):=e^{-E_{\{x\}}^{\Phi }(\gamma +\varepsilon _{x})}\rho
(x),\;x\in X,\gamma \in \Gamma  \label{5eq55}
\end{equation}
It was shown in \cite[Theorem 5.3]{AR90} (in the case $X=\QTR{mathbb}{R}^{d}$%
) that the components of the Dirichlet form $(\mathcal{E}_{\sigma _{\gamma
}}^{X},\mathcal{D}^{\sigma _{\gamma }})$ corresponding to the measure $%
\sigma _{\gamma }$ are closable on $L^{2}(\QTR{mathbb}{R}^{d},\sigma
_{\gamma })$ if and only if $\sigma _{\gamma }$ is absolutely continuous
with respect to Lebesgue measure on $\QTR{mathbb}{R}^{d}$ and the
Radon-Nikodym derivative satisfies some condition, see (\ref{5eq51}) below
for details. This result allows us to prove the closability of $(\mathcal{E}%
_{\mu }^{\Gamma },\mathcal{F}C_{b}^{\infty }(\mathcal{D},\Gamma ))$ on $%
L^{2}(\mu )$. Let us first recall the above mentioned result.

\begin{theorem}
\label{5eq56}(cf.~\cite[Theorem 5.3]{AR90}) Let $\nu $ by a probability
measure on $(\QTR{mathbb}{R}^{d},\mathcal{B}(\QTR{mathbb}{R}^{d})$, $d\in 
\QTR{mathbb}{N}$ and let $\mathcal{D}^{\nu }$ denote the $\nu $-classes
determined by $\mathcal{D}$. Then the forms $(\mathcal{E}_{\nu ,i},\mathcal{D%
}^{\nu })$ defined by 
\[
\mathcal{E}_{\nu ,i}(u,v):=\int_{\QTR{mathbb}{R}^{d}}\frac{\partial u}{%
\partial x_{i}}\frac{\partial v}{\partial x_{i}}d\nu ,\;u,v\in \mathcal{D},
\]
are well-defined and closable on $L^{2}(\QTR{mathbb}{R}^{d},\nu )$ for $%
1\leq i\leq d$ if and only if $\nu $ is absolutely continuous with respect
to Lebesgue measure $\lambda ^{d}$ on $\QTR{mathbb}{R}^{d}$, and the
Radon-Nikodym derivative $\rho =d\nu /d\lambda ^{d}$ satisfies the
condition: 
\begin{eqnarray}
&&\mathrm{for\;any\;}1\leq i\leq d\;\mathrm{and\;}\lambda ^{d-1}\mathrm{%
-a.e.\;}x\in \left\{ y\in \QTR{mathbb}{R}^{d-1}|\int_{\QTR{mathbb}{R}}\rho
_{y}^{(i)}(s)d\lambda ^{1}(s)>0\right\} ,  \nonumber \\
&&\rho _{x}^{(i)}=0\;\lambda ^{1}\mathrm{-a.e.\;on\;}\QTR{mathbb}{%
R\backslash R(\rho }_{x}^{(i)})\mathrm{,\;where\;}\QTR{mathbb}{\rho }%
_{x}^{(i)}(s):=\rho (x_{1},\ldots ,x_{i-1},s,x_{i},\ldots ,x_{d}),  \nonumber
\\
&&s\in \QTR{mathbb}{R}\mathrm{,\;if\;}x=(x_{1},\ldots ,x_{d-1})\in 
\QTR{mathbb}{R}^{d-1}\mathrm{,\;and\;where}  \label{5eq51}
\end{eqnarray}
\begin{equation}
R(\rho _{x}^{(i)}):=\left\{ t\in \QTR{mathbb}{R|}\int_{t-\varepsilon
}^{t+\varepsilon }\frac{1}{\rho _{x}^{(i)}(s)}ds<\infty \;\mathrm{for\;some\;%
}\varepsilon >0\right\} .  \label{5eq52}
\end{equation}
\end{theorem}

There is an obvious generalization of Theorem \ref{5eq56} to the case where
a Riemannian manifold $X$ is replacing $\QTR{mathbb}{R}^{d}$, to be
formulated in terms of local charts. Since here we are only interested in
the ``if part'' of Theorem \ref{5eq56}, we now recall a slightly weaker
sufficient condition for closability in the general case where $X$ is a
manifold as before.

\begin{theorem}
\label{5eq58}Suppose $\sigma _{1}=\rho _{1}\cdot m$, where $\rho
_{1}:X\rightarrow \QTR{mathbb}{R}_{+}$ is $\mathcal{B}(X)$-measurable such
that 
\begin{equation}
\rho _{1}=0\;m\mathrm{-a.e.~on\;}X\backslash \!\left\{ x\in
X|\!\int_{\Lambda _{x}}\frac{1}{\rho _{1}}dm<\infty \mathrm{%
\;for\;some\;open\;neighbourhood\;}\Lambda _{x}\mathrm{\;of\;}x\right\} .
\label{5eq57}
\end{equation}
Then $(\mathcal{E}_{\sigma _{1}}^{X},\mathcal{D}^{\sigma _{1}})$ defined by 
\[
\mathcal{E}_{\sigma _{1}}^{X}(u,v):=\int_{X}\langle \nabla ^{X}u(x),\nabla
^{X}v(x)\rangle _{T_{x}X}\,d\sigma _{1}(x);\;u,v\in \mathcal{D},
\]
is closable on $L^{2}(\sigma _{1}).$
\end{theorem}

The proof is a straightforward adaptation of the line of arguments in 
\cite[Chap.~II, Subsection 2a]{MR92} (see also \cite[Theorem 4.2]{ABR89} for
details). We emphasize that (\ref{5eq57}) e.g.~always holds, if $\rho _{1}$
is lower semicontinuous, and that neither $\nu $ in Theorem \ref{5eq56} nore 
$\sigma _{1}$ in Theorem \ref{5eq58} is required to have full support, so
e.g.~$\rho _{1}$ is not necessarily strictly positive $m$-a.e.~on $X$.

We are now ready to prove the closability of $(\mathcal{E}_{\mu }^{\Gamma },%
\mathcal{F}C_{b}^{\infty }(\mathcal{D},\Gamma ))$ on $L^{2}(\mu )$ under the
above assumption.

\begin{theorem}
\label{5eq65}Let $\mu \in \mathcal{G}_{gc}^{1}(\sigma ,\Phi )$. Suppose that
for $\mu $-a.e.~$\gamma \in \Gamma $ the function $\rho _{\gamma }$ defined
in (\ref{5eq55}) satisfies (\ref{5eq57}) (resp.~(\ref{5eq51}) in case $X=%
\QTR{mathbb}{R}^{d}$). Then the form $(\mathcal{E}_{\mu }^{\Gamma },\mathcal{%
F}C_{b}^{\infty }(\mathcal{D},\Gamma ))$ is closable on $L^{2}(\mu )$.
\end{theorem}

\noindent \textbf{Proof.} Let $(F_{n})_{n\in \QTR{mathbb}{N}}$ be a sequence
in $\mathcal{F}C_{b}^{\infty }(\mathcal{D},\Gamma )$ such that $%
F_{n}\rightarrow 0$, $n\rightarrow \infty $ in $L^{2}(\mu )$ and 
\begin{equation}
\mathcal{E}_{\mu }^{\Gamma }(F_{n}-F_{m},F_{n}-F_{m})\stackunder{%
n,m\rightarrow \infty }{\longrightarrow }0.  \label{5eq54}
\end{equation}
We have to show that 
\begin{equation}
\mathcal{E}_{\mu }^{\Gamma }(F_{n_{k}},F_{n_{k}})\stackunder{k\rightarrow
\infty }{\longrightarrow }0  \label{5eq59}
\end{equation}
for some subsequence $(n_{k})_{k\in \QTR{mathbb}{N}}$. Let $(n_{k})_{k\in 
\QTR{mathbb}{N}}$ be a subsequence such that 
\[
\left( \int_{\Gamma }F_{n_{k}}^{2}d\mu \right) ^{1/2}+\mathcal{E}_{\mu
}^{\Gamma }(F_{n_{k+1}}-F_{n_{k}},F_{n_{k+1}}-F_{n_{k}})^{1/2}<\frac{1}{2^{k}%
}\mathrm{\;for\;all\;}k\in \QTR{mathbb}{N}. 
\]
Then 
\begin{eqnarray*}
\infty &>&\sum_{k=1}^{\infty }\mathcal{E}_{\mu }^{\Gamma
}(F_{n_{k+1}}-F_{n_{k}},F_{n_{k+1}}-F_{n_{k}})^{1/2} \\
&\geq &\sum_{k=1}^{\infty }\int_{\Gamma }\left( \int_{X}|\nabla ^{X}\nabla
^{P}(F_{n_{k+1}}-F_{n_{k}})(x,\gamma )|_{T_{x}X}^{2}\,e^{-E_{\{x\}}^{\Phi
}(\gamma +\varepsilon _{x})}d\sigma (x)\right) ^{1/2}d\mu (\gamma ) \\
&=&\int_{\Gamma }\sum_{k=1}^{\infty }\left( \int_{X}|\nabla ^{X}\nabla
^{P}(F_{n_{k+1}}-F_{n_{k}})(x,\gamma )|_{T_{x}X}^{2}\,\rho _{\gamma
}(x)dm(x)\right) ^{1/2}d\mu (\gamma ),
\end{eqnarray*}
where we used Theorem \ref{5eq2} and (\ref{5eq55}). From the last expression
we obtain that 
\begin{eqnarray}
\!\! &&\!\!\sum_{k=1}^{\infty }\mathcal{E}_{\sigma _{\gamma
}}^{X}(u_{n_{k+1}}^{(\gamma )}-u_{n_{k}}^{(\gamma )},u_{n_{k+1}}^{(\gamma
)}-u_{n_{k}}^{(\gamma )})^{1/2}  \label{5eq60} \\
\!\! &=&\!\!\sum_{k=1}^{\infty }\left( \int_{X}\!|\nabla ^{X}\nabla
^{P}(F_{n_{k+1}}-F_{n_{k}})(x,\gamma )|_{T_{x}X}^{2}\rho _{\gamma
}(x)dm(x)\right) ^{1/2}\!\!\!<\!\infty \mathrm{\;for\;}\mu \mathrm{-a.e.\;}%
\gamma \in \Gamma ,  \nonumber
\end{eqnarray}
where for $k\in \QTR{mathbb}{N}$, $\gamma \in \Gamma $, 
\[
u_{n_{k}}^{(\gamma )}(x):=F_{n_{k}}(\gamma +\varepsilon
_{x})-F_{n_{k}}(\gamma ),\;x\in X. 
\]
Note that $u_{n_{k}}^{(\gamma )}\in \mathcal{D}$. (\ref{5eq60}) implies that
for $\mu $-a.e.~$\gamma \in \Gamma $%
\begin{equation}
\mathcal{E}_{\sigma _{\gamma }}^{X}(u_{n_{k}}^{(\gamma )}-u_{n_{l}}^{(\gamma
)},u_{n_{k}}^{(\gamma )}-u_{n_{l}}^{(\gamma )})\stackunder{k,l\rightarrow
\infty }{\longrightarrow }0.  \label{5eq61}
\end{equation}
Let $\Lambda \subset \mathcal{O}_{c}(X)$.\bigskip

\noindent \textbf{Claim 1}: For $\mu $-a.e.~$\gamma \in \Gamma $%
\[
\int_{X}(u_{n_{k}}^{(\gamma )}(x))^{2}1\!\!1_{\Lambda }(x)d\sigma _{\gamma
}(x)\stackunder{k\rightarrow \infty }{\longrightarrow }0. 
\]
To prove Claim 1 we first note that for $\mu $-a.e.~$\gamma \in \Gamma $%
\[
\sigma _{\gamma }(\Lambda )<\infty , 
\]
as follows immediately from (\ref{5eq4}) (taking $h(\gamma
,x):=1\!\!1_{\Lambda }(x)$ for $x\in X$, $\gamma \in \Gamma $), since $\mu
\in \mathcal{G}_{gc}^{1}(\sigma ,\Phi )$. Therefore, for $\mu $-a.e.~$\gamma
\in \Gamma $%
\begin{equation}
\int_{X}F_{n_{k}}^{2}(\gamma )1\!\!1_{\Lambda }(x)d\sigma _{\gamma
}(x)=F_{n_{k}}^{2}(\gamma )\sigma _{\gamma }(\Lambda )\stackunder{%
k\rightarrow \infty }{\longrightarrow }0.  \label{5eq62}
\end{equation}
Furthermore, by (\ref{5eq4}) 
\begin{eqnarray*}
&&\int_{\Gamma }\int_{X}F_{n_{k}}^{2}(\gamma +\varepsilon
_{x})1\!\!1_{\Lambda }(x)d\sigma _{\gamma }(x)(1+\gamma (\Lambda ))^{-1}d\mu
(\gamma ) \\
&=&\int_{\Gamma }F_{n_{k}}^{2}(\gamma )\int_{X}\frac{1\!\!1_{\Lambda }(x)}{%
1+\gamma (\Lambda )-1\!\!1_{\Lambda }(x)}\gamma (dx)d\mu (\gamma ) \\
&\leq &\int_{\Gamma }F_{n_{k}}^{2}(\gamma )d\mu (\gamma )<\frac{1}{2^{k}},
\end{eqnarray*}
because the integral w.r.t.$~\gamma $ is dominated by 1 for all $\gamma \in
\Gamma $. Hence 
\begin{eqnarray*}
\infty &>&\sum_{k=1}^{\infty }\left( \int_{\Gamma
}\int_{X}F_{n_{k}}^{2}(\gamma +\varepsilon _{x})1\!\!1_{\Lambda }(x)d\sigma
_{\gamma }(x)(1+\gamma (\Lambda ))^{-1}d\mu (\gamma )\right) ^{1/2} \\
&\geq &\int_{\Gamma }\sum_{k=1}^{\infty }\left( \int_{X}F_{n_{k}}^{2}(\gamma
+\varepsilon _{x})1\!\!1_{\Lambda }(x)d\sigma _{\gamma }(x)\right)
^{1/2}(1+\gamma (\Lambda ))^{-1}d\mu (\gamma ).
\end{eqnarray*}
Therefore, for $\mu $-a.e.~$\gamma \in \Gamma $%
\begin{equation}
\int_{X}F_{n_{k}}^{2}(\gamma +\varepsilon _{x})1\!\!1_{\Lambda }(x)d\sigma
_{\gamma }(x)\stackunder{k\rightarrow \infty }{\longrightarrow }0.
\label{5eq63}
\end{equation}
Then Claim 1 follows by (\ref{5eq62}) and (\ref{5eq63}).\bigskip

\noindent \textbf{Claim 2}: For $\mu $-a.e.~$\gamma \in \Gamma $%
\[
|\nabla ^{X}u_{n_{k}}^{(\gamma )}|_{TX}\stackunder{k\rightarrow \infty }{%
\longrightarrow }0\;\;\;\sigma _{\gamma }\mathrm{-a.e.} 
\]
To prove Claim 2 we first note that clearly (\ref{5eq60}) implies that for $%
\mu $-a.e.~$\gamma \in \Gamma $%
\begin{equation}
\mathcal{E}_{1\!\!1_{\Lambda }\sigma _{\gamma }}^{X}(u_{n_{k}}^{(\gamma
)}-u_{n_{l}}^{(\gamma )},u_{n_{k}}^{(\gamma )}-u_{n_{l}}^{(\gamma )})%
\stackunder{k,l\rightarrow \infty }{\longrightarrow }0.  \label{5eq64}
\end{equation}
Hence we can apply Theorem \ref{5eq58} (resp.~\ref{5eq56}) to $\rho
_{1}:=1\!\!1_{\Lambda }\rho _{\gamma }$ and conclude by Claim 1 and (\ref
{5eq64}) that for $\mu $-a.e.~$\gamma \in \Gamma $%
\[
\mathcal{E}_{1\!\!1_{\Lambda }\sigma _{\gamma }}^{X}(u_{n_{k}}^{(\gamma
)},u_{n_{k}}^{(\gamma )})\stackunder{k\rightarrow \infty }{\longrightarrow }%
0, 
\]
hence by (\ref{5eq60}) 
\[
1\!\!1_{\Lambda }|\nabla ^{X}u_{n_{k}}^{(\gamma )}|_{TX}\stackunder{%
k\rightarrow \infty }{\longrightarrow }0\;\;\;\sigma _{\gamma }\mathrm{-a.e.}
\]
Since $\Lambda $ was arbitrary, Claim 2 is proven.

From Claim 2 we now easily deduce (\ref{5eq59}) by (\ref{5eq3}) and Fatou's
Lemma as follows: 
\begin{eqnarray*}
\mathcal{E}_{\mu }^{\Gamma }(F_{n_{k}},F_{n_{k}}) &\leq &\int_{\Gamma }%
\stackunder{l\rightarrow \infty }{\lim \inf }\int_{X}|\nabla
^{X}(u_{n_{k}}^{(\gamma )}-u_{n_{l}}^{(\gamma )})|_{TX}^{2}\,d\sigma
_{\gamma }(x)d\mu (\gamma ) \\
&\leq &\stackunder{l\rightarrow \infty }{\lim \inf }\mathcal{E}_{\mu
}^{\Gamma }(F_{n_{k}}-F_{n_{l}},F_{n_{k}}-F_{n_{l}}),
\end{eqnarray*}
which by (\ref{5eq54}) can be made arbitrarily small for $k$ large enough.$%
\hfill \blacksquare $

\begin{remark}
The above method to prove closability of pre-Dirichlet forms on
configuration spaces $\Gamma _{X}$ extends immediately to the case where $X$
is replaced by an infinite dimensional ``manifold'' such as the loop space
(cf.~\cite{MR97}).
\end{remark}

\begin{example}
\label{5eq66}Let $X=\QTR{mathbb}{R}^{d}$ with the Euclidean metric and $%
\sigma :=z\cdot m$, $z\in (0,\infty )$. A pair potential is a $\mathcal{B}(%
\QTR{mathbb}{R}^{d})$-measurable function $\phi :\QTR{mathbb}{R}%
^{d}\rightarrow \QTR{mathbb}{R}\cup \{\infty \}$ such that $\phi (-x)=\phi
(x)$. Any pair potential $\phi $ defines a potential $\Phi =\Phi _{\phi }$
in the sense of Section~\ref{5eq6} as follows: we set $\Phi (\gamma ):=0$, $%
|\gamma |\neq 2$ and $\Phi (\gamma ):=\phi (x-y)$ for $\gamma
=\{x,y\}\subset \QTR{mathbb}{R}^{d}$. For such pair potentials $\phi $ the
condition in Theorem~\ref{5eq65} ensuring closability of $(\mathcal{E}_{\mu
}^{\Gamma },\mathcal{F}C_{b}^{\infty }(\mathcal{D},\Gamma ))$ on $L^{2}(\mu )
$ for $\mu \in \mathcal{G}_{gc}^{1}(\sigma ,\Phi )$ can be now easily
formulated as follows: for $\mu $-a.e.~$\gamma \in \Gamma $ and $m$-a.e.~$%
x\in \{y\in \QTR{mathbb}{R}^{d}|\sum_{y^{\prime }\in \gamma \backslash
\{y\}}|\phi (y-y^{\prime })|<\infty \}$ it holds that $\int_{V_{x}}e^{%
\sum_{y^{\prime }\in \gamma \backslash \{y\}}\phi (y-y^{\prime
})}m(dy)<\infty $ for some open neighbourhood $V_{x}$ of $x$. This condition
trivially holds e.g.~if \textrm{supp}$\phi $ is compact, $\{\phi <\infty \}$
is open, and $\phi ^{+}\in L_{loc}^{\infty })(\{\phi <\infty \};m)$. If even 
$\mu \in \mathcal{G}_{gc}^{t}(z,\phi )$ and $\phi $ satisfies the
assumptions in Proposition~7.1, then it suffices to merely assume that $%
\{\phi <\infty \}$ is open and $\phi ^{+}\in L_{loc}^{\infty }(\{\phi
<\infty \};m)$. This follows by an elementary consideration.
\end{example}

\begin{remark}
\label{5eq67}

\begin{enumerate}
\item  \label{5eq68}We emphasize that Example~\ref{5eq66} generalizes the
closability result in \cite{O96}, though an a-priori bigger domain for $%
\mathcal{E}_{\mu }^{\Gamma }$ is considered there. However, Theorems~\ref
{5eq56}-\ref{5eq65} are also valid for this bigger domain. The proofs are
exactly the same.

\item  \label{5eq69}We also like to emphasize that similarly to Example~\ref
{5eq66} one proves the closability of $(\mathcal{E}_{\mu }^{\Gamma },%
\mathcal{F}C_{b}^{\infty }(\mathcal{D},\Gamma ))$ (or with a larger domain
in \cite{O96}) on $L^{2}(\mu )$ for $\mu \in \mathcal{G}_{gc}^{1}(\sigma
,\Phi )$ in the case of multi-body potentials $\phi $.\bigskip 
\end{enumerate}
\end{remark}

\noindent \textbf{Acknowledgments}\medskip

We would like to thank Tobias Kuna for helpful discussions. Financial
support of the INTAS-Project 378, PRAXIS Programme through CITMA, Funchal,
and TMR Nr.~ERB4001GT957046 are gratefully acknowledged.

\end{document}